\documentclass[11pt]{article}
\usepackage[T2A]{fontenc}
\usepackage[cp1251]{inputenc}
\usepackage{amsmath}
\usepackage{amssymb}
\def\a{\mathbf{a}}
\def\b{\mathbf{b}}
\def\betta{\boldsymbol{\beta}}
\def\p{\mathbf{p}}
\def\m{\mathbf{m}}
\def\f{\mathbf{f}}

\def\N{\mathbb{N}}
\def\Z{\mathbb{Z}}

\def\C{\mathbb{C}}

\def\SS{\mathbf{S}}
\def\etta{\boldsymbol{\eta}}
\def\zetta{\boldsymbol{\zeta}}
\def\llambda{\boldsymbol{\lambda}}
\def\ggamma{\boldsymbol{\gamma}}

\newtheorem{lemma}{\hspace*{\parindent}Lemma}
\newtheorem{theorem}{\hspace*{\parindent}Theorem}
\newtheorem{corollary}{\hspace*{\parindent}Corollary}
\title{Alternative approach to Miller-Paris transformations and their extensions}
\author{D.B.\:Karp$^{\rm a,b}$\footnote{Corresponding author. E-mail: D. Karp -- \emph{dimkrp@gmail.com},
E.\:Prilepkina --  \emph{pril-elena@yandex.ru}}~~and
E.G.\:Prilepkina$^{\rm a,b}$
\\[10pt]
\\
\small{\textit{$\phantom{1}^a$Far Eastern Federal University, 8
Sukhanova street, Vladivostok, 690950, Russia}}
\\
\small{\textit{$\phantom{1}^b$Institute of Applied Mathematics,
FEBRAS, 7 Radio Street, Vladivostok,  690041, Russia}}}
\date{}

\begin{document}
\maketitle

\begin{center}
\parbox{12cm}{
\small\textbf{Abstract.}
Miller-Paris transformations are extensions of Euler's transformations for the Gauss hypergeometric functions to generalized hypergeometric functions of higher-order having integral parameter differences (IPD).  In our recent work we computed the degenerate versions of these transformations corresponding to the case when one parameter difference is equal to a negative integer.  The purpose of this paper is to present an independent new derivation of both the general and the degenerate forms of Miller-Paris transformations. In doing so we employ the generalized Stieltjes transform representation of the generalized hypergeometric functions and some partial fraction expansions. Our approach leads to different forms of the characteristic polynomials, one of them appears noticeably simpler than the original form due to Miller and Paris. We further present two extensions of the degenerate transformations to the generalized hypergeometric functions with additional free parameters and additional parameters with negative integral differences.}
\end{center}

\bigskip


\bigskip


\bigskip

\section{Introduction and preliminaries}

Transformation, reduction and summation formulas for hypergeometric functions is a vast subject with rich history dating back to Leonard Euler.
Among  important applications of such formulas (let alone hypergeometric functions in general) are quantum physics  \cite{RJRJR}, non-equilibrium statistical
physics \cite{KMT2018} and many other fields \cite{Seaborn}. The main developments up to the end of 20th century can be found, for instance, in the books \cite{AAR,BealsWong,LukeBook}.  One particular example useful for simplifying sums that arise in theoretical physics (such as Racah coefficients) is the following
summation formula established by Minton \cite{Minton} in 1970:
\begin{equation}\label{eq:Minton}
{}_{r+2}F_{r+1}\!\!\left(\left.\begin{matrix}-k,b,f_1+m_1,\ldots,f_r+m_r\\b+1,f_1,\ldots,f_r\end{matrix}\right|1\right)=\frac{k!}{(b+1)_{k}}\frac{(f_1-b)_{m_1}\cdots(f_r-b)_{m_r}}{(f_1)_{m_1}\cdots(f_r)_{m_r}}, ~k\ge{m},~k\in\N,
\end{equation}
and slightly generalized by Karlsson \cite{Karlson} who replaced $-k$ by an arbitrary complex number $a$ satisfying $\Re(1-a-m)>0$ to get
 \begin{equation}\label{eq:Karl1}
{}_{r+2}F_{r+1}\!\!\left(\left.\begin{matrix}a,b,f_1+m_1,\ldots,f_r+m_r\\b+1,f_1,\ldots,f_r\end{matrix}\right|1\right)
=\frac{\Gamma(b+1)\Gamma(1-a)}{\Gamma(b+1-a)}\frac{(f_1-b)_{m_1}\cdots(f_r-b)_{m_r}}{(f_1)_{m_1}\cdots(f_r)_{m_r}}.
\end{equation}
Here and throughout the paper ${_{p}F_q}$ stands for the generalized hypergeometric function (see \cite[Section~2.1]{AAR}, \cite[Section~5.1]{LukeBook}, \cite[Sections 16.2-16.12]{NIST} or \cite[Chapter~12]{BealsWong}), $(a)_k=\Gamma(a+k)/\Gamma(a)$ is rising factorial and $\Gamma(z)$ is Euler's gamma function.  These formulas attracted attention to generalized hypergeometric function with integral parameter differences, for which Michael Schlosser subsequently introduced the acronym IPD, motivated by the title of Karlsson's paper \cite{Karlson}.  These summation formulas were generalized and extended in many directions: Gasper \cite{Gasper} deduced a $q$-analogue and a generalization of Minton's and Karlsson's formulas; Chu \cite{Chu1,Chu2} found extensions to bilateral hypergeometric and $q$-hypergeometric series; their results were re-derived by simpler means and further generalized by Schlosser \cite{Schlosser2}, who also found multidimensional extensions  to hypergeometric functions associated with root systems \cite{Schlosser1}.  For further developments in this directions, see also \cite{Rosengren1,Rosengren2}.   We also mention an interesting work \cite{LVW} by Letessier, Valent  and Wimp, where an order reduction for the differential equation satisfied by the generalized hypergeometric functions with some integral parameter differences was established.

Another surge in interest to IPD-type hypergeometric functions is related with transformation formulas for such functions, generalizing the classical Euler's transformations for the Gauss function ${}_2F_{1}$ and Kummer's transformation for the confluent hypergeometric functions ${}_1F_{1}$.  Unlike Minton-Karlsson formulas dealing the generalized hypergeometric functions  evaluated at $1$ these transformations are certain identities for these functions evaluated at an arbitrary value of the argument.  They were developed in a series  of papers published over last 15 years, the most general form was presented in a seminal paper \cite{MP2013} by  Miller and Paris. For a vector of positive integers $\m=(m_1,\ldots,m_r)$, $m=m_1+m_2+\ldots+m_r$, and a complex vector $\f=(f_1,\ldots,f_r)$ these transformations are given by  \cite[Theorem~1]{KRP2014}
\begin{equation}\label{eq:KRPTh1-1}
{}_{r+2}F_{r+1}\left.\!\!\left(\!\begin{matrix}a, b,\f+\m\\c,\f\end{matrix}\right\vert x\right)
=(1-x)^{-a}{}_{m+2}F_{m+1}\left.\!\!\left(\!\begin{matrix}a, c-b-m, \zetta+1\\c, \zetta\end{matrix}\right\vert\frac{x}{x-1}\right)
\end{equation}
if $(c-b-m)_m\ne0$, and, if also $(c-a-m)_m\ne0$ and $(1+a+b-c)_m\ne0$, then
\begin{equation}\label{eq:KRPTh1-2}
{}_{r+2}F_{r+1}\left.\!\!\left(\!\begin{matrix}a, b,\f+\m\\c,\f\end{matrix}\right\vert x\right)
=(1-x)^{c-a-b-m}{}_{m+2}F_{m+1}\left.\!\!\left(\!\begin{matrix}c-a-m, c-b-m, \etta+1\\c, \etta\end{matrix}\right\vert x\right).
\end{equation}
Here  the vector $\zetta=\zetta(c,b,\f)=(\zeta_1,\ldots,\zeta_m)$ comprises the roots of the polynomial
\begin{equation}\label{eq:Qm}
Q_m(t)=Q(b,c,\f,\m;t)=\frac{1}{(c-b-m)_{m}}\sum\limits_{k=0}^{m}(b)_kC_{k,r}(t)_{k}(c-b-m-t)_{m-k},
\end{equation}
where $C_{0,r}=1$, $C_{m,r}=1/(\f)_{\m}$, $(\f)_{\m}=(f_1)_{m_1}\cdots(f_r)_{m_r}$, and
\begin{equation}\label{eq:Ckr}
C_{k,r}=C_{k,r}(\f,\m)=\frac{1}{(\f)_{\m}}\sum\limits_{j=k}^m\sigma_j\mathbf{S}_j^{(k)}=\frac{(-1)^k}{k!}{}_{r+1}F_{r}\!\left(\begin{matrix}-k,\f+\m\\\f\end{matrix}\right).
\end{equation}
In this formula and below we routinely omit the argument $1$ from the generalized hypergeometric function: ${_{p}F_q}(\a;\b):={_{p}F_q}(\a;\b;1)$.
The numbers $\sigma_j$ ($0\leq{j}\leq{m}$) are defined via the generating function
$$
(f_1+x)_{m_1}\ldots(f_r+x)_{m_r}=\sum\limits_{j=0}^m\sigma_jx^j,
$$
and $\mathbf{S}_j^{(k)}$ stands for the Stirling's number of the second kind.   A simple rearrangement of Pochhammer's symbols leads to an alternative form of the polynomial $Q_m(t)$
as given in \cite[(3.7)]{KPITSF2018}:
\begin{equation}\label{eq:sumMil2}
Q(b,c,\f,\m;t)=\frac{(c-b-t-m)_m}{(c-b-m)_m}\sum\limits_{k=0}^m
{}_{r+1}F_{r}\left.\!\!\left(\begin{matrix}-k,\f+\m\\\f\end{matrix}\right.\right)
\frac{(t)_k(b)_k}{(1+t+b-c)_k k!}.
\end{equation}
Further, $\etta=(\eta_1,\ldots,\eta_m)$ in (\ref{eq:KRPTh1-2}) are the roots of
\begin{equation}\label{eq:Qmhat}
\hat{Q}_m(t)=\sum\limits_{k=0}^{m}\frac{(-1)^kC_{k,r}(a)_k(b)_k(t)_k}{(c-a-m)_k(c-b-m)_k}
{}_{3}F_{2}\!\left(\begin{matrix}-m+k,t+k,c-a-b-m\\c-a-m+k,c-b-m+k\end{matrix}\right).
\end{equation}
See \cite{KRP2014,Miller2005,MP2011,MP2012R,MS2010} and references therein for further details.

Both formulas (\ref{eq:KRPTh1-1}) and (\ref{eq:KRPTh1-2}) fail when $c=b+p$, $p\in\{1,\ldots,m\}$.  In our recent paper \cite{KPReluts2018} we used a careful limit transition to derive the degenerate forms of the transformations (\ref{eq:KRPTh1-1}) and (\ref{eq:KRPTh1-2}) for such values of $c$.  When evaluated at $x=1$ these degenerate transformations lead to extensions of Minton-Karlsson summation formulas (\ref{eq:Minton}), (\ref{eq:Karl1}) which we also investigated in \cite{KPITSF2018} using several different techniques.  The main purpose of this paper is to present an alternative derivation of both  the general transformations (\ref{eq:KRPTh1-1}), (\ref{eq:KRPTh1-2})  and  their degenerate forms found in \cite{KPReluts2018}.  Our derivation of the general case is presented in Section~2 and is based on the representation of the generalized hypergeometric function by the generalized Stieltjes transform of a particular type of Meijer's $G$ function, namely $G_{p,p}^{p,0}$, which for historical reasons we prefer to call the Meijer-N{\o}rlund function. Details regarding this representation,  its history and numerous applications can be found in \cite{KarpJMS2015,KLJAT2017,KPITSF2017}.  Our approach leads to different forms of the characteristic polynomials (\ref{eq:Qm}), (\ref{eq:Qmhat}).  Comparison with those new forms yields an identity for finite hypergeometric sums which may be difficult to obtain directly. The degenerate forms and their extensions are presented in Section~3. Their derivation hinges on certain simple partial fraction decompositions.
The results differ from those in \cite{KPReluts2018}:  here the negative parameter difference $-p$ (recall that $b-c=-p$) may take any (negative) integer value regardless of whether degeneration happens or not in the corresponding general Miller-Paris transformation. We further present two extensions: to several negative parameters differences instead of one and to a pair of additional unrestricted parameters on top and bottom of the generalized hypergeometric function.

\section{Miller-Paris transformations: general case}
Before proceeding to the main results let us introduce some notation. Let $\N$ and $\C$ denote the natural and complex numbers, respectively; further, put
$$
\Gamma(\a)=\Gamma(a_1)\Gamma(a_2)\cdots\Gamma(a_p),~~(\a)_n=(a_1)_n(a_2)_n\cdots(a_p)_n,~~\a+\mu=(a_1+\mu,a_2+\mu,\dots,a_p+\mu).
$$
Inequalities like $\Re(\a)>0$ and properties like $-\a\notin\N_0:=\N\cup\{0\}$ will be understood element-wise. The symbol $\a_{[k]}$ will stand for the vector $\a$ with omitted $k$-th component. The function $G^{m,n}_{p,q}$ is Meijer's $G$ function (see \cite[section~5.2]{LukeBook}, \cite[16.17]{NIST}, \cite[8.2]{PBM3} or \cite[Chapter~12]{BealsWong}).

We begin with a lemma expressing the Meijer-N{\o}rlund function $G^{p,0}_{p,p}$ with integral parameter differences in terms of beta density times a rational function.
\begin{lemma}\label{th:G-finitesum}
Let $\m=(m_1,\ldots,m_r)\in\N^r$, $m=m_1+m_2+\ldots+m_r$, $\f=(f_1,\ldots,f_r)\in\C^r$ and $b,c\in\C$.
Then
\begin{equation}\label{eq:G-finitesum}
G_{r+1,r+1}^{r+1,0}\!\left(\!t~\vline\begin{array}{l}c,\f\\b,\f+\m\end{array}\!\!\right)
=\frac{t^b(1-t)^{c-b-1}}{\Gamma(c-b)}
\sum\limits_{k=0}^m D_k(c-b-k)_{k}\frac{t^k}{(t-1)^{k}},
\end{equation}
where
\begin{equation}\label{eq:Dk}
D_k=D_k(\f,\m,b)=\sum\limits_{j=k}^m\alpha_j\SS_{j}^{(k)}
=\frac{(-1)^k(\f-b)_{\m}}{k!}
{}_{r+1}F_{r}\!\left(\begin{matrix}-k,1-\f+b\\1-\f+b-\m\end{matrix}\right).
\end{equation}
The numbers $\alpha_j$ are defined via the generating function
\begin{equation}\label{eq:alpha}
(\f-b-t)_{\m}=\sum\limits_{j=0}^{m}\alpha_jt^j,
\end{equation}
and $ \SS_{j}^{(k)}$ stands for the Stirling number of the second kind.
\end{lemma}

\textbf{Proof.}  By the well-known expansion of the Meijer-N{\o}rlund function $G_{r+1,r+1}^{r+1,0}$, see, for instance \cite[(2.4)]{KPSIGMA},
if $f_i-f_j\notin\Z,$ $f_i-c\notin\Z$ we have
\begin{multline*}
G_{r+1,r+1}^{r+1,0}\!\left(\!t~\vline\begin{array}{l}c,\f\\b,\f+\m\end{array}\!\!\right)=\frac{t^b\Gamma(\f+\m-b)}{\Gamma(c-b)\Gamma(\f-b)}
{}_{r+1}F_{r}\left.\!\!\left(\begin{matrix}1-c+b,1-\f+b\\1-\f-\m+b\end{matrix}\right\vert t\right)
\\
=\frac{t^b\Gamma(\f+\m-b)}{\Gamma(c-b)\Gamma(\f-b)}\sum\limits_{n=0}^{\infty}\frac{(1-c+b)_n(1-\f+b)_n}{(1-\f-\m+b)_nn!}t^n.
\end{multline*}
Next, using $\Gamma(\f-b)(\f-b)_{\m}=\Gamma(\f+\m-b)$ and
$$
\frac{(1-\f+b)_n}{(1-\f-\m+b)_n}=\frac{(\f-b-n)_{\m}}{(\f-b)_{\m}}=\frac{1}{(\f-b)_{\m}}\sum\limits_{k=0}^{m}\alpha_kn^k,
$$
where $\alpha_k$ is defined in (\ref{eq:alpha}), we obtain:
\begin{multline*}
G_{r+1,r+1}^{r+1,0}\!\left(\!t~\vline\begin{array}{l}c,\f\\b,\f+\m\end{array}\!\!\right)
\\
=\frac{t^b\Gamma(\f+\m-b)}{\Gamma(c-b)\Gamma(\f-b)(\f-b)_{\m}}\sum\limits_{n=0}^{\infty}\frac{(1-c+b)_n}{n!}t^n\sum\limits_{k=0}^{m}\alpha_kn^k
=\frac{t^b}{\Gamma(c-b)}\sum\limits_{k=0}^{m}\alpha_k\sum\limits_{n=0}^{\infty}\frac{(1-c+b)_n}{n!}t^nn^k
\\
=\frac{t^b}{\Gamma(c-b)}\sum\limits_{k=0}^m\alpha_k\sum\limits_{l=0}^k \SS_{k}^{(l)}\frac{(1-c+b)_lt^l}{(1-t)^{1-c+b+l}}
=\frac{t^b(1-t)^{c-b-1}}{\Gamma(c-b)}\sum\limits_{l=0}^m\frac{(-1)^l(c-b-l)_{l}t^l}{(1-t)^{l}}\sum\limits_{k=l}^m\alpha_k\SS_{k}^{(l)}.
\end{multline*}
To get the pre-ultimate  equality we applied the definition of the Stirling numbers of the second kind via
 $n^k=\sum\nolimits_{l=0}^k\SS_{k}^{(l)}[n]_l,$ $[n]_l=n(n-1)\ldots(n-l+1)$ and the next relation (with $\delta=1-c+b$):
\begin{multline*}
\sum\limits_{n=0}^{\infty} \frac{(\delta)_n t^n}{n!}n^k=\sum\limits_{n=0}^{\infty} \frac{(\delta)_n t^n}{n!}\sum\limits_{l=0}^k\SS_{k}^{(l)}[n]_l=
\sum\limits_{l=0}^k\SS_{k}^{(l)}\sum\limits_{n=l}^{\infty} \frac{(\delta)_n t^n}{(n-l)!}=\\
\sum\limits_{l=0}^k\SS_{k}^{(l)}\sum\limits_{n=0}^{\infty} \frac{(\delta)_{n+l} t^{n+l}}{n!}=
\sum\limits_{l=0}^k\SS_{k}^{(l)}\sum\limits_{n=0}^{\infty} \frac{(\delta)_l(\delta+l)_n t^{n+l}}{n!}=\sum\limits_{l=0}^k\SS_{k}^{(l)}\frac{(\delta)_lt^l}{(1-t)^{\delta+l}}.
\end{multline*}
Thus, we have proved formula (\ref{eq:G-finitesum}) with the expression for $D_k$ given by the first equality in (\ref{eq:Dk}). It remains to prove the second equality in (\ref{eq:Dk}). To this end we will borrow the technique from \cite[Theorem~2]{MP2012R}. By Taylor's theorem:
$$
R(t)=(\f-b-t)_{\m}=\sum\limits_{j=0}^{m}\alpha_jt^j=\sum\limits_{j=0}^{m}R^{(j)}(0)\frac{t^j}{j!}.
$$
On the other hand, from the theory of finite differences:
$$
\Delta^{k}R(t)=\sum\limits_{j=0}^{k}(-1)^{k-j}\binom{k}{j}R(t+j)=k!\sum\limits_{j=k}^{m}\SS_{j}^{(k)}\frac{R^{(j)}(t)}{j!}.
$$
Comparing with the first formula in (\ref{eq:Dk}) we get:
$$
D_k=\sum\limits_{j=k}^{m}\SS_{j}^{(k)}\frac{R^{(j)}(0)}{j!}=\frac{1}{k!}\Delta^{k}R(0)
=\frac{1}{k!}\sum\limits_{j=0}^{k}(-1)^{k-j}\binom{k}{j}R(j)
=\frac{1}{k!}\sum\limits_{j=0}^{k}(-1)^{k-j}\binom{k}{j}(\f-b-j)_{\m}.
$$
Substituting the identity
$$
(\f-b-j)_{\m}=(\f-b)_{\m}\frac{(1-\f+b)_{j}}{(1-f+b-\m)_{j}}
$$
into the above expression after simple rearrangement and in view of the formula
$$
\frac{1}{(k-j)!}=(-1)^j\frac{(-k)_j}{k!}.
$$
we get the second equality in (\ref{eq:Dk}).$\hfill\square$

As a corollary we get an expansion of ${}_{r+2}F_{r+1}$ with $r$ positive integer parameter differences.

\begin{corollary}\label{cr:GenHyp2F1exp}
The following expansion holds:
\begin{equation}\label{eq:GenHyp2F1exp}
{}_{r+2}F_{r+1}\left.\!\!\left(\begin{matrix}a,b,\f+\m\\c,\f\end{matrix}\right\vert x\right)
=\frac{1}{(\f)_{\m}}\sum\limits_{k=0}^m(-1)^kD_k(b)_k{}_{2}F_{1}\left.\!\!\left(\begin{matrix}a,b+k\\c\end{matrix}\right\vert x\right)
\end{equation}
\end{corollary}
\textbf{Proof.} Indeed, by \cite[(2)]{KLJAT2017} and (\ref{eq:G-finitesum})
\begin{multline*}
\frac{\Gamma(\f+\m)}{\Gamma(\f)}{}_{r+2}F_{r+1}\left.\!\!\left(\begin{matrix}a,b,\f+\m\\c,\f\end{matrix}\right\vert x\right)
=\frac{\Gamma(c)}{\Gamma(b)}\int\limits_0^{1}G_{r+1,r+1}^{r+1,0}\!\left(\!t~\vline\begin{array}{l}c,\f\\b,\f+\m\end{array}\!\!\right)\frac{dt}{t(1-xt)^a}
\\
=\frac{\Gamma(c)}{\Gamma(b)\Gamma(c-b)}
\sum\limits_{l=0}^m\alpha_l\sum\limits_{k=0}^l \SS_{l}^{(k)}(1-c+b)_k\int\limits_{0}^{1}\frac{t^{b+k-1}(1-t)^{c-b-k-1}}{(1-xt)^a}dt
\\
=\frac{\Gamma(c)}{\Gamma(b)\Gamma(c-b)}
\sum\limits_{l=0}^m\alpha_l\sum\limits_{k=0}^l \SS_{l}^{(k)}(1-c+b)_k
\frac{\Gamma(b+k)\Gamma(c-b-k)}{\Gamma(c)}{}_{2}F_{1}\left.\!\!\left(\begin{matrix}a,b+k\\c\end{matrix}\right\vert x\right)
\\
=\sum\limits_{l=0}^m\alpha_l\sum\limits_{k=0}^l \SS_{l}^{(k)}
(-1)^k(b)_k{}_{2}F_{1}\left.\!\!\left(\begin{matrix}a,b+k\\c\end{matrix}\right\vert x\right)
=\sum\limits_{k=0}^m(-1)^k(b)_k{}_{2}F_{1}\left.\!\!\left(\begin{matrix}a,b+k\\c\end{matrix}\right\vert x\right)\sum\limits_{l=k}^m\SS_{l}^{(k)}\alpha_l.~\square
\end{multline*}

Using Corollary~\ref{cr:GenHyp2F1exp} we can easily recover the first Miller-Paris transformation (\ref{eq:KRPTh1-1}) (see \cite[(1.3)]{MP2013}, \cite[Theorem~1]{KRP2014})  in the following form.
\begin{theorem}\label{th:MP1}
Let $\m=(m_1,\ldots,m_r)\in\N^r$, $m=m_1+m_2+\ldots+m_r$, $\f=(f_1,\ldots,f_r)\in\C^r$ and $a,b,c\in\C$ be such that $(c-b-m)_{m}\ne0$.
Then for all $x\in\C\!\setminus\![1,\infty)$ we have:
\begin{equation}\label{eq:MP1}
{}_{r+2}F_{r+1}\left.\!\!\left(\begin{matrix}a,b,\f+\m\\c,\f\end{matrix}\right\vert x\right)
=(1-x)^{-a}{}_{m+2}F_{m+1}\left.\!\!\left(\begin{matrix}a,c-b-m,\zetta+1\\c,\zetta\end{matrix}\right\vert \frac{x}{x-1}\right),
\end{equation}
where $\zetta=(\zeta_1,\ldots,\zeta_m)$ are the roots of the polynomial
$$
P_m(x)=\frac{1}{(c-b-m)_m}\sum\limits_{k=0}^m(b)_{k}(1-c+b)_{k}D_{k}(c-b-m-x)_{m-k}
$$
with $D_k$ defined in \emph{(\ref{eq:Dk})}.
\end{theorem}
\textbf{Proof}.  Apply Euler-Pfaff transformation  \cite[(2.2.6)]{AAR} to the Gauss function ${}_{2}F_{1}$ on the right hand side of (\ref{eq:GenHyp2F1exp}) and calculate:
\begin{multline*}
\frac{\Gamma(\f+\m)}{\Gamma(\f)}{}_{r+2}F_{r+1}\left.\!\!\left(\begin{matrix}a,b,\f+\m\\c,\f\end{matrix}\right\vert x\right)
=(1-x)^{-a}\sum\limits_{l=0}^m(-1)^lD_{l}(b)_{l}{}_{2}F_{1}\left.\!\!\left(\begin{matrix}a,c-b-l\\c\end{matrix}\right\vert \frac{x}{x-1}\right)
\\
=(1-x)^{-a}\sum\limits_{l=0}^m(-1)^lD_{l}(b)_{l}\sum\limits_{j=0}^{\infty}\frac{(a)_j(c-b-l)_j}{(c)_jj!}\frac{x^j}{(x-1)^j}
\\
=(1-x)^{-a}\sum\limits_{j=0}^{\infty}\frac{(a)_{j}x^j}{(c)_{j}j!(x-1)^j}\sum\limits_{l=0}^m(-1)^lD_{l}(b)_{l}(c-b-l)_j
\\
=(1-x)^{-a}\sum\limits_{j=0}^{\infty}\frac{(a)_{j}x^j}{(c)_{j}j!(x-1)^j}\sum\limits_{l=0}^m(-1)^lD_{l}(b)_{l}\frac{(c-b-m)_{j}(c-b-m+j)_{m-l}}{(c-b-m)_{m-l}}
\\
=(1-x)^{-a}\sum\limits_{j=0}^{\infty}\frac{(a)_{j}(c-b-m)_{j}x^j}{(c)_{j}j!(x-1)^j}\tilde{P}_m(j)
\\
=(1-x)^{-a}\sum\limits_{j=0}^{\infty}\frac{(a)_{j}(c-b-m)_{j}x^j}{(c)_{j}j!(x-1)^j}A_m(j-\zeta_1)\cdots(j-\zeta_m)
\\
=(1-x)^{-a}\sum\limits_{j=0}^{\infty}\frac{(a)_{j}(c-b-m)_{j}x^j}{(c)_{j}j!(x-1)^j}
A_m(-\zeta_1)(-\zeta_2)\cdots(-\zeta_m)\frac{(1-\zeta_1)_j\cdots(1-\zeta_m)_j}{(-\zeta_1)_j\cdots(-\zeta_m)_j},
\\
=(1-x)^{-a}\tilde{P}_m(0){}_{m+2}F_{m+1}\left.\!\!\left(\begin{matrix}a,c-b-m,-\zetta+1\\c,-\zetta\end{matrix}\right\vert \frac{x}{x-1}\right),
\end{multline*}
where
$$
\tilde{P}_m(x)=\sum\limits_{l=0}^m(-1)^lD_{l}(b)_{l}\frac{(c-b-m+x)_{m-l}}{(c-b-m)_{m-l}}
=A_m(x-\zeta_1)\cdots(x-\zeta_m),
$$
and the identity
$$
(c-b-l)_{j}=\frac{(c-b-m)_{j}(c-b-m+j)_{m-l}}{(c-b-m)_{m-l}}
$$
has been used.  It remains to define $P_m(x)=\tilde{P}_m(-x)$ and note that $\tilde{P}_m(0)$ must equal $(\f)_{\m}$ by taking $x=0$ in the resulting identity. $\hfill\square$

By comparing the Miller-Paris formula (\ref{eq:KRPTh1-1}) with (\ref{eq:MP1}) it is clear that $P_m(x)$ must be a constant multiple of $Q_m(x)$ defined in (\ref{eq:Qm}).

We give a direct proof of this fact below.

\begin{lemma}\label{lm:PmQm}
We have
$$
P_m(x)=(\f)_{\m}Q(b,c,\f,\m;x).
$$
\end{lemma}
\textbf{Proof.}  A combination of  \cite[Theorem~2.3]{KPITSF2018} with \cite[Theorem~3.2]{KPITSF2018} gives the following identity for $Q_m$:
$$
\frac{(\f-b)_{\m}(1-c+x)_{m}}{(\f)_{\m}(1-c+b)_{m}}Q(1-c+b,1-x+b,1-\f+b-\m,\m;b)=Q(b,c,\f,\m;x).
$$
Using this identity and  definition (\ref{eq:Qm}) we obtain
\begin{multline*}
P_m(x)=\frac{1}{(c-b-m)_m}\sum\limits_{k=0}^m(b)_{k}(1-c+b)_{k}D_{k}(c-b-m-x)_{m-k}
\\
=\frac{(\f-b)_{\m}}{(c-b-m)_m}\sum\limits_{k=0}^m(b)_{k}(1-c+b)_{k}
\frac{(-1)^k}{k!}{}_{r+1}F_{r}\!\left(\begin{matrix}-k,1-\f+b\\1-\f+b-\m\end{matrix}\right)
(c-b-m-x)_{m-k}
\\
=\frac{(\f-b)_{\m}(1-c+x)_m}{(1-c+b)_m}Q(1-c+b,1-x+b,1-\f+b-\m,\m;b)=(\f)_{\m}Q(b,c,\f,\m;x).~~~~\square
\end{multline*}

Our version of the second Miller-Paris transformation (\ref{eq:KRPTh1-2}) (see \cite[Theorem~4]{MP2013},\cite[Theorem~1]{KRP2014}) is the following theorem.

\begin{theorem}\label{th:MP2}
Let $\m=(m_1,\ldots,m_r)\in\N^r$, $m=m_1+m_2+\ldots+m_r$, $\f=(f_1,\ldots,f_r)\in\C^r$ and $a,b,c\in\C$ be such that $(c-a-m)_{m}\ne0$, $(c-b-m)_{m}\ne0$, $(1+a+b-c)_m\ne0$.
Then for all $x\in\C\!\setminus\![1,\infty)$ we have:
\begin{equation}\label{eq:MP2}
{}_{r+2}F_{r+1}\left.\!\!\left(\begin{matrix}a,b,\f+\m\\c,\f\end{matrix}\right\vert x\right)
=(1-x)^{c-a-b-m}{}_{m+2}F_{m+1}\left.\!\!\left(\begin{matrix}c-a-m,c-b-m,\etta+1\\c,\etta\end{matrix}\right\vert x\right),
\end{equation}
where $\etta=(\eta_1,\ldots,\eta_m)$ are the roots of the polynomial
\begin{equation}\label{eq:hatP}
\hat{P}_m(t)=\sum\limits_{k=0}^m\frac{(-1)^k(a)_k(-b-m)_k(t)_k(c-a-m-t)_{m-k}}{(c-a-m)_m(c-b-m)_kk!}
{}_{r+2}F_{r+1}\!\left(\begin{matrix}-k,b,\f+\m\\b+m-k+1,\f\end{matrix}\right).
\end{equation}
\end{theorem}

The proof of this theorem will require the following two lemmas which might be of independent interest.

\begin{lemma}\label{lm:sumlemma}
For any nonnegative integers $0\le{i}\le{k}\le{m}$ the following summation formula holds\emph{:}
\begin{equation}\label{eq:sumlemma}
\sum\limits_{j=i}^{k}(-k)_j(\alpha-m+j)_{m-i}\frac{(-j)_i}{j!}=(-1)^i(-k)_i\frac{(-m)_{k}(\alpha-m)_{m}}{(-m)_{i}(\alpha-m)_{k}}.
\end{equation}
\end{lemma}
\textbf{Proof.} Writing $S$ for the left hand side of (\ref{eq:sumlemma}) and using the straightforward identities
$$
(-j)_i=\frac{(-1)^ij!}{(j-i)!},~~~~(\alpha-m+j)_{m-i}=\frac{(\alpha-i)_{j}(\alpha-m)_{m-i}}{(\alpha-m)_j},~~~~(\beta)_{k+r}=(\beta)_k(\beta+k)_{r},
$$
we compute by changing the summation index to $n=j-i$:
\begin{multline*}
S=(-1)^i(\alpha-m)_{m-i}\sum\limits_{j=i}^{k}\frac{(-k)_j(\alpha-i)_{j}}{(\alpha-m)_j(j-i)!}
=(-1)^i(\alpha-m)_{m-i}\sum\limits_{n=0}^{k-i}\frac{(-k)_{n+i}(\alpha-i)_{n+i}}{(\alpha-m)_{n+i}n!}
\\
=(-1)^i(\alpha-m)_{m-i}\sum\limits_{n=0}^{k-i}\frac{(-k)_{i}(-k+i)_{n}(\alpha-i)_{i}(\alpha)_{n}}{(\alpha-m)_{i}(\alpha-m+i)_{n}n!}
\\
=\frac{(-1)^i(-k)_{i}(\alpha-i)_{i}(\alpha-m)_{m-i}}{(\alpha-m)_{i}}\sum\limits_{n=0}^{k-i}\frac{(-k+i)_{n}(\alpha)_{n}}{(\alpha-m+i)_{n}n!}
\\
=\frac{(-1)^i(-k)_{i}(\alpha-i)_{i}(\alpha-m)_{m-i}}{(\alpha-m)_{i}}\frac{(-m+i)_{k-i}}{(\alpha-m+i)_{k-i}}
=(-1)^i(-k)_i(i-m)_{k-i}\frac{(\alpha-m)_{m}}{(\alpha-m)_{k}}.
\end{multline*}
The pre-ultimate equality here is the celebrated Chu-Vandermonde identity  \cite[Corollary~2.2.3]{AAR}. It remains to apply $(i-m)_{k-i}=(-m)_{k}/(-m)_{i}$. $\hfill\square$

\begin{lemma}\label{lm:hypersum}
For any nonnegative integers $0\le{k}\le{m}$ the following summation formula holds\emph{:}
\begin{equation}\label{eq:hypersum}
\sum\limits_{i=0}^{k}\frac{(-k)_i(b)_i}{(-m)_{i}i!}{}_{r+1}F_{r}\!\left(\begin{matrix}-i,\f+\m\\\f\end{matrix}\right)
=\frac{(-b-m)_k}{(-m)_k}{}_{r+2}F_{r+1}\!\left(\begin{matrix}-k,b,\f+\m\\b+m-k+1,\f\end{matrix}\right).
\end{equation}
\end{lemma}
\textbf{Proof.} Let $\epsilon$ be a small positive number. According to (\ref{eq:sumMil2}) with $t=-k$, $c=b+m-k+1+\epsilon$,
$$
\sum\limits_{i=0}^{k}\frac{(-k)_i(b)_i}{(-m-\epsilon)_{i}i!}{}_{r+1}F_{r}\!\left(\begin{matrix}-i,\f+\m\\\f\end{matrix}\right)
=\frac{(1-k+\epsilon)_m}{(1+\epsilon)_m}Q_m(b;b+m-k+1+\epsilon;\f;\m;-k)
$$
Further, by \cite[Theorem~3.2]{KPITSF2018},
\begin{multline*}
\frac{(1-k+\epsilon)_m}{(1+\epsilon)_m}Q_m(b;b+m-k+1+\epsilon;\f;\m;-k)
\\
=\frac{(1-k+\epsilon)_m}{(1+\epsilon)_m}\frac{\Gamma(b+m+1+\epsilon)\Gamma(1-k+\epsilon)}{\Gamma(b+m-k+1+\epsilon)\Gamma(1+\epsilon)}
{}_{r+2}F_{r+1}\!\left(\begin{matrix}-k,b,\f+\m\\b+m-k+1+\epsilon,\f\end{matrix}\right)
\\
=\frac{\Gamma(1-k+\epsilon+m)}{\Gamma(1+\epsilon+m)}\frac{\Gamma(b+m+1+\epsilon)}{\Gamma(b+m-k+1+\epsilon)}
{}_{r+2}F_{r+1}\!\left(\begin{matrix}-k,b,\f+\m\\b+m-k+1+\epsilon,\f\end{matrix}\right).
\end{multline*}
Letting $\epsilon\to0$ and using $m!/(m-k)!=(-1)^k(-m)_k$ and $\Gamma(b+m+1)/\Gamma(b+m-k+1)=(-1)^k(-b-m)_k$ we arrive at (\ref{eq:hypersum}).$\hfill\square$

\medskip

\noindent\textbf{Remark.} This lemma could also be derived from \cite[Theorem~4]{MP2012R} but as conditions of this theorem are violated here, some special treatment would still be needed. We also prefer to give an independent proof based entirely on our results.

\textbf{Proof of Theorem~\ref{th:MP2}.}  Applying (\ref{eq:KRPTh1-1}) with $a$ and $b$ interchanged to the right hand side of the same formula we immediately get (\ref{eq:MP2}) with the characteristic polynomial given by
$$
\hat{P}_m(t)=\frac{1}{(c-a-m)_m}\sum\limits_{k=0}^m(a)_k(c-a-m-t)_{m-k}(t)_k\frac{(-1)^k}{k!}{}_{r+1}F_{r}\!\left(\begin{matrix}-k,\zetta+1\\\zetta\end{matrix}\right),
$$
where $\zetta=\zetta(b,c,\f)=(\zeta_1,\ldots,\zeta_m)$ are the roots of the polynomial $Q_m(b,c,\f;t)$ defined in (\ref{eq:Qm}).
Next, we calculate using $Q_m(0)=1$ and (\ref{eq:sumlemma}):
\begin{multline*}
{}_{r+1}F_{r}\!\left(\begin{matrix}-k,\zetta+1\\\zetta\end{matrix}\right)=1+\frac{(-k)(\zetta+1)}{(\zetta)1!}
+\cdots+\frac{(-k)_{k}(\zetta+1)_{k}}{(\zetta)_{k}k!}
\\
1+\frac{(-k)(\zetta+1)}{(\zetta)1!}+\cdots+\frac{(-k)_{k}(\zetta+k)}{(\zetta)k!}
=1+\frac{(-k)Q_m(-1)}{Q_m(0)1!}+\cdots+\frac{(-k)_{k}Q_m(-k)}{Q_m(0)k!}
\\
=\frac{1}{(c-b-m)_{m}}\sum\limits_{j=0}^{k}\frac{(-k)_j}{j!}\sum\limits_{i=0}^{m}(b)_i(-j)_{i}(c-b-m+j)_{m-i}\frac{(-1)^i}{i!}
{}_{r+1}F_{r}\!\left(\begin{matrix}-i,\f+\m\\\f\end{matrix}\right)
\\
=\frac{1}{(c-b-m)_{m}}\sum\limits_{i=0}^{k}\frac{(-1)^i(b)_i}{i!}{}_{r+1}F_{r}\!\left(\begin{matrix}-i,\f+\m\\\f\end{matrix}\right)\sum\limits_{j=i}^{k}\frac{(-k)_j}{j!}(-j)_{i}(c-b-m+j)_{m-i}
\\
=\frac{(-m)_{k}}{(c-b-m)_{k}}\sum\limits_{i=0}^{k}\frac{(-k)_i(b)_i}{(-m)_{i}i!}{}_{r+1}F_{r}\!\left(\begin{matrix}-i,\f+\m\\\f\end{matrix}\right)
=\frac{(-b-m)_k}{(c-b-m)_{k}}{}_{r+2}F_{r+1}\!\left(\begin{matrix}-k,b,\f+\m\\b+m-k+1,\f\end{matrix}\right),
\end{multline*}
where we applied (\ref{eq:hypersum}) in the last equality.  Substituting this formula into the above expression for $\hat{P}_m$ we arrive at (\ref{eq:hatP}). $\hfill\square$

\begin{corollary}
The following identity is true: $\hat{Q}_m(t)=\hat{P}_m(t)$,  where $\hat{Q}_m$ is defined in \emph{(\ref{eq:Qmhat})} and $\hat{P}_m$ is defined in \emph{(\ref{eq:hatP})}.
\end{corollary}

\textbf{Proof.} Indeed, comparing (\ref{eq:KRPTh1-2}) and (\ref{eq:MP2}) we see that $\hat{Q}_m(t)$ and $\hat{P}_m(t)$ have the same roots and thus may only differ by a nonzero multiplicative constant. However, it is straightforward that $\hat{Q}_m(0)=\hat{P}_m(0)=1$ and the claim follows. $\hfill\square$

Note that the identity $\hat{Q}_m(t)=\hat{P}_m(t)$ represents a non-trivial hypergeometric transformation which seems to be hard to obtain directly.

\section{Miller-Paris transformations: degenerate case}

As we mentioned in the introduction and the statements of Theorems~\ref{th:MP1} and \ref{th:MP2}, formulas (\ref{eq:MP1}) and (\ref{eq:MP2}) fail when $c=b+p$, $p\in\{1,\ldots,m\}$.  The purpose of this section is to present two transformations valid when  $c=b+p$ with arbitrary $p\in\N$.  Hence, they cover both degenerate and non-degenerate cases.  Some of the coefficients appearing in these transformations can be expressed in terms of N{\o}rlund's coefficients $g_n(\a;\b)$ which were introduced by N{\o}rlund in \cite[(1.33)]{Norlund} and investigated in our papers \cite[section~2.2]{KPSIGMA}, \cite[Property~6]{KLJAT2017} and \cite[section~2]{KPITSF2018}.  For completeness we also give a short and slightly different
account here.  The  functions $g_n(\a;\b)$, $n\in\N_0$, are polynomials separately symmetric in the components of the vectors $\a=(a_1,\ldots,a_{p-1})$ and $\b=(b_1,\ldots,b_p)$.
They can be defined either via the power series generating function \cite[(1.33)]{Norlund}, \cite[(11)]{KLJAT2017}
\begin{equation}\label{eq:Norlund}
G^{p,0}_{p,p}\!\left(\!1-z~\vline\begin{array}{l}\b\\\a,0\end{array}\!\!\right)=\frac{z^{\nu_p-1}}{\Gamma(\nu_p)}
\sum\limits_{n=0}^{\infty}\frac{g_n(\a;\b)}{(\nu_p)_n}z^n,
\end{equation}
where $\nu_p=\nu_p(\a;\b)=\sum_{j=1}^{p}b_j-\sum_{j=1}^{p-1}a_j$, or via the inverse factorial generating function \cite[(2.21)]{Norlund}
$$
\frac{\Gamma(z+\nu_p)\Gamma(z+\a)}{\Gamma(z+\b)}=\sum\limits_{n=0}^{\infty}\frac{g_n(\a;\b)}{(z+\nu_p)_n}.
$$
As, clearly, $\nu_p(\a+\alpha;\b+\alpha)=\nu_p(\a;\b)+\alpha$, we have (by changing $z\to{z+\alpha}$)
$$
\frac{\Gamma(z+\alpha+\nu_p)\Gamma(z+\alpha+\a)}{\Gamma(z+\alpha+\b)}=\sum\limits_{n=0}^{\infty}\frac{g_n(\a+\alpha;\b+\alpha)}{(z+\nu_p+\alpha)_n}
=\sum\limits_{n=0}^{\infty}\frac{g_n(\a;\b)}{(z+\alpha+\nu_p)_n}.
$$
Hence, $g_n(\a+\alpha;\b+\alpha)=g_n(\a;\b)$ for any $\alpha$. N{\o}rlund found two different recurrence relations for $g_n(\a;\b)$ (one in $p$ and one in $n$).
The simplest of them reads \cite[(2.7)]{Norlund}
\begin{equation}\label{eq:Norlundcoeff}
g_n(\a,\alpha;\b,\beta)=\sum\limits_{s=0}^{n}\frac{(\beta-\alpha)_{n-s}}{(n-s)!}(\nu_p-\alpha+s)_{n-s}g_s(\a;\b),~~~p=1,2,\ldots,
\end{equation}
with the  initial values $g_0(-;b_1)=1$, $g_n(-;b_1)=0$, $n\ge1$.  This recurrence was solved by N{\o}rlund  \cite[(2.11)]{Norlund} as follows:
\begin{equation}\label{eq:Norlund-explicit}
g_n(\a;\b)=\sum\limits_{0\leq{j_{1}}\leq{j_{2}}\leq\cdots\leq{j_{p-2}}\leq{n}}
\prod\limits_{m=1}^{p-1}\frac{(\psi_m+j_{m-1})_{j_{m}-j_{m-1}}}{(j_{m}-j_{m-1})!}(b_{m+1}-a_{m})_{j_{m}-j_{m-1}},
\end{equation}
where $\psi_m=\sum_{i=1}^{m}(b_i-a_i)$, $j_0=0$, $j_{p-1}=n$.   Another recurrence relation for $g_n(\a;\b)$ discovered by N{\o}rlund \cite[(1.28)]{Norlund} has order $p$ in the variable $n$ and coefficients polynomial in $n$. Details can be found in \cite[section~2.2]{KPSIGMA}.  The first three coefficients  are given by \cite[Theorem~2]{KPSIGMA}:
$$
g_0(\a;\b)=1,~~~~g_1(\a;\b)=\sum_{m=1}^{p-1}(b_{m+1}-a_m)\psi_m,
$$
$$
g_2(\a;\b)=\frac{1}{2}\sum_{m=1}^{p-1}(b_{m+1}-a_m)_2(\psi_m)_2+\sum_{k=2}^{p-1}(b_{k+1}-a_{k})(\psi_{k}+1)\sum_{m=1}^{k-1}(b_{m+1}-a_m)\psi_m.
$$
For $p=2$ and $p=3$ and arbitrary $n\in\N_0$ explicit expressions for $g_n(\a;\b)$ discovered  by N{\o}rlund \cite[eq.(2.10)]{Norlund} are:
\begin{equation}\label{eq:gnp2p3}
\begin{split}
&g_n(a;\b)=\frac{(b_1-a)_{n}(b_2-a)_{n}}{n!}~~\text{for}~p=2;
\\
&g_n(\a;\b)=\frac{(\nu_3-b_2)_n(\nu_3-b_3)_n}{n!}
{}_3F_2\left(\!\!\begin{array}{l}-n,b_1-a_1,b_1-a_2\\\nu_3-b_2,\nu_3-b_3\end{array}\!\!\right)~~\text{for}~p=3,
\end{split}
\end{equation}
where $\nu_m:=\sum_{j=1}^{m}b_j-\sum_{j=1}^{m-1}a_j$.  The right hand side here is invariant with respect to the permutation of the elements of $\b$. Finally, for $p=4$ we have \cite[p.12]{KPSIGMA}
$$
g_n(\a;\b)=\frac{(\nu_4-b_3)_n(\nu_4-b_4)_n}{n!}
\sum\limits_{k=0}^{n}\frac{(-n)_k(\nu_2-a_{2})_k(\nu_2-a_{3})_k}{(\nu_4-b_3)_k(\nu_4-b_4)_k}
{}_{3}F_{2}\!\left(\!\begin{array}{l}-k,b_1-a_{1},b_2-a_{1}\\\nu_2-a_{2}, \nu_2-a_{3}\end{array}\!\!\right).
$$

Next, define $W_{m-1}(n)=W_{m-1}(b,\f,\m;n)=\sum\nolimits_{k=0}^{m-1}\delta_kn^k$ to be the polynomial of degree $m-1$ given by
\begin{equation}\label{eq:pol}
W_{m-1}(n)=W_{m-1}(b,\f,\m;n)=\left(\frac{(\f+\m)_n}{(\f)_n}-\frac{(\f-b)_\m}{(\f)_\m}\right)\frac{(b)_n}{(b+1)_n}= \frac{b((\f+n)_\m-(\f-b)_\m)}{(b+n)(\f)_\m}.
\end{equation}
The following theorem gives two extensions of the Karlsson's formula (\ref{eq:Karl1}) to arbitrary argument.
\begin{theorem} \label{thm:Fpred}
The following transformation formulas hold:
\begin{equation}\label{eq:Kar111}
(1-x)^a{}_{r+2}F_{r+1}\left.\!\!\!\left(\begin{matrix}a, b,\f+\m\\ b+1,\f\end{matrix}\right\vert x\right)
=\frac{(\f-b)_\m}{(\f)_\m}{}_{2}F_{1}\left.\!\!\!\left(\begin{matrix}1,a\\b+1\end{matrix}\right\vert\frac{x}{x-1}\right)
+\sum\limits_{l=0}^{m-1}Y_l\frac{(a)_lx^l}{(1-x)^l}
\end{equation}
and
\begin{equation}\label{eq:Kar112}
(1-x)^{a-1}{}_{r+2}F_{r+1}\left.\!\!\!\left(\begin{matrix}a, b,\f+\m\\ b+1,\f\end{matrix}\right\vert x\right)
=\frac{(\f-b)_\m}{(\f)_\m}{}_{2}F_{1}\left.\!\!\!\left(\begin{matrix}1, b+1-a\\b+1\end{matrix}\right\vert x\right)+\sum\limits_{l=0}^{m-1}Y_l\frac{(a)_lx^l}{(1-x)^{l+1}}.
\end{equation}
Here
\begin{multline}\label{eq:Bl}
Y_l=Y_l(b,\f,\m)=\sum\limits_{k=l}^{m-1}\delta_k\SS_{k}^{(l)}
=\frac{(-1)^{l}}{l!}{}_{r+2}F_{r+1}\!\!\left(\begin{matrix}-l, b,\f+\m\\ b+1,\f\end{matrix}\right)
-\frac{(-1)^{l}(\f-b)_\m}{(b+1)_l(\f)_\m}
\\
=\frac{(-1)^{m-l-1}b}{(\f)_\m}\sum\limits_{i=0}^{m-1-l}(-1)^ig_{m-1-l-i}(-\f;-\f-\m,l)(1-b)_i,
\end{multline}
where $\delta_k$ are the coefficients of the polynomial $W_{m-1}(x)$ defined in \emph{(\ref{eq:pol})} and $g_{n}(\cdot;\cdot)$ are N{\o}rlund's coefficients given in \emph{(\ref{eq:Norlund-explicit})}.
\end{theorem}

\textbf{Proof.}  First, we show that for $|x|<1$ the following equality holds:
\begin{equation}\label{eq:KarNew1}
{}_{r+2}F_{r+1}\left.\!\!\!\left(\begin{matrix}a, b,
\f+\m\\ b+1,\f\end{matrix}\right\vert x\right)
=\frac{(\f-b)_\m}{(\f)_\m}{}_{2}F_{1}\left.\!\!\!\left(\begin{matrix}a, b\\ b+1\end{matrix}\right\vert x\right)+\sum\limits_{l=0}^{m-1}\frac{(a)_lx^l}{(1-x)^{a+l}}\sum\limits_{k=l}^{m-1}\delta_k\SS_{k}^{(l)}.
\end{equation}
Indeed, in view of (\ref{eq:pol}),
\begin{multline*}
{}_{r+2}F_{r+1}\left.\!\!\!\left(\begin{matrix}a, b,
\f+\m\\ b+1,\f\end{matrix}\right\vert x\right)=\frac{(\f-b)_\m}{(\f)_\m}\sum\limits_{n=0}^\infty \frac{(a)_n (b)_n x^n}{n!(b+1)_n}+\sum\limits_{n=0}^\infty \frac{(a)_n  x^n W_{m-1}(n)}{n!}
\\
=\frac{(\f-b)_\m}{(\f)_\m}{}_{2}F_{1}\left.\!\!\!\left(\begin{matrix}a, b
\\ b+1\end{matrix}\right\vert x\right)+\sum\limits_{k=0}^{m-1}\delta_k \sum\limits_{n=0}^{\infty} \frac{(a)_n x^n}{n!}n^k.
\end{multline*}
Using the definition of the Stirling numbers $n^k=\sum\nolimits_{l=0}^k\SS_{k}^{(l)}[n]_l$ in terms of falling factorials $[n]_l=n(n-1)\ldots(n-l+1)$, we get
\begin{multline*}
\sum\limits_{n=0}^{\infty} \frac{(a)_n x^n}{n!}n^k=\sum\limits_{n=0}^{\infty} \frac{(a)_n x^n}{n!}\sum\limits_{l=0}^k\SS_{k}^{(l)}[n]_l
=\sum\limits_{l=0}^k\SS_{k}^{(l)}\sum\limits_{n=l}^{\infty} \frac{(a)_n x^n}{(n-l)!}
\\
=\sum\limits_{l=0}^k\SS_{k}^{(l)}\sum\limits_{n=0}^{\infty} \frac{(a)_{n+l} x^{n+l}}{n!}=
\sum\limits_{l=0}^k\SS_{k}^{(l)}\sum\limits_{n=0}^{\infty} \frac{(a)_l(a+l)_n x^{n+l}}{n!}=\sum\limits_{l=0}^k\SS_{k}^{(l)}(a)_l\frac{x^l}{(1-x)^{a+l}},
\end{multline*}
which implies (\ref{eq:KarNew1}) after exchanging the order of summations.  It remains to apply Euler's transformations to the ${}_2F_1$ on the right hand side of (\ref{eq:KarNew1}) to get (\ref{eq:Kar111}) and (\ref{eq:Kar112}) with $Y_l$ given by the first formula in (\ref{eq:Bl}).

To obtain the second expression for $Y_l$ recall that
$$
W_{m-1}(x)=\frac{b[(\f+x)_{\m}-(\f-b)_{\m}]}{(b+x)(\f)_\m}=\sum\limits_{k=0}^{m-1}\frac{W_{m-1}^{(k)}(0)}{k!}x^k,
$$
yielding $\delta_k=W_{m-1}^{(k)}(0)/k!$.  Using the technique from \cite[Theorem~2]{MP2012R} we now apply the following formula
from the theory of finite differences:
$$
\Delta^{l}W_{m-1}(x)=\sum\limits_{j=0}^{l}(-1)^{l-j}\binom{l}{j}W_{m-1}(x+j)=l!\sum\limits_{j=l}^{m-1}\SS_{j}^{(k)}\frac{W_{m-1}^{(j)}(x)}{j!}.
$$
In view of  \eqref{eq:pol} we have
\begin{multline*}
\frac{1}{l!}\Delta^{l}W_{m-1}(0)=\sum\limits_{k=l}^{m-1}\delta_k\SS_{k}^{(l)}
=\frac{1}{l!}\sum\limits_{j=0}^{l}(-1)^{l-j}\binom{l}{j}W_{m-1}(j)
\\
=\frac{1}{l!}\sum\limits_{j=0}^{l}\frac{(-1)^{l-j}(b)_j}{(b+1)_j}\binom{l}{j}
\left[\frac{(\f+\m)_j}{(\f)_j}-\frac{(\f-b)_\m}{(\f)_\m}\right]
=\frac{(-1)^{l}}{l!}\sum\limits_{j=0}^{l}\frac{(-l)_j(b)_j}{(b+1)_jj!}\left[\frac{(\f+\m)_j}{(\f)_j}-\frac{(\f-b)_\m}{(\f)_\m}\right]
\\
=\frac{(-1)^{l}}{l!}{}_{r+2}F_{r+1}\!\!\left(\begin{matrix}-l, b,\f+\m\\ b+1,\f\end{matrix}\right)
-\frac{(-1)^{l}}{l!}\frac{(\f-b)_\m}{(\f)_\m}{}_{2}F_{1}\!\!\left(\begin{matrix}-l, b\\ b+1\end{matrix}\right)
\\
=\frac{(-1)^{l}}{l!}{}_{r+2}F_{r+1}\!\!\left(\begin{matrix}-l, b,\f+\m\\ b+1,\f\end{matrix}\right)
-\frac{(-1)^{l}(\f-b)_\m}{(b+1)_l(\f)_\m},
\end{multline*}
where we employed the relation
$$
\binom{l}{j}=(-1)^j\frac{(-l)_j}{j!}
$$
in the second line and the Chu-Vandermonde identity in the last equality.  Further, according to \cite[Theorem~2.1]{KPITSF2018}
$$
{}_{r+2}F_{r+1}\left.\!\!\left(\begin{matrix}-l,b,\f+\m\\b+1,\f\end{matrix}\right.\right)
=\frac{l!}{(b+1)_{l}}\frac{(\f-b)_\m}{(\f)_\m}-\frac{(-1)^ml!b}{(\f)_\m}q_l,
$$
where $q_l=\sum_{i=0}^{m-l-1}g_{m-l-i-1}(b-\f;b-\f-\m,b+l)(b-i)_i$, and N{\o}rlund's coefficient $g_n(\cdot;\cdot)$ is defined in (\ref{eq:Norlund-explicit}).
Substituting and using the shifting property $g_n(\a+\alpha;\b+\alpha)=g_n(\a;\b)$  of N{\o}rlund's coefficients, we get:
$$
\sum\limits_{k=l}^{m-1}\delta_k\SS_{k}^{(l)}=\frac{(-1)^{m-l-1}b}{(\f)_\m}\sum\limits_{i=0}^{m-l-1}g_{m-l-i-1}(-\f;-\f-\m,l)(b-i)_i
$$
which is equivalent to the second formula in (\ref{eq:Bl}).$\hfill\square$

\medskip

\textbf{Remark.}
If $\Re(1-a-m)>0$ an application of the Gauss summation formula to ${}_2F_1$ on the right hand side of (\ref{eq:KarNew1}) results in Karlsson's formula (\ref{eq:Karl1}).

\medskip

Theorem~\ref{thm:Fpred} can be  generalized as follows. Suppose $\b=(b_1,\ldots,b_l)$ is a complex vector, $\p=(p_1,\ldots,p_l)$ is a vector of positive integers, $p=p_1+p_2+\ldots+p_l$, and all elements of the vector $\betta=(b_1,b_1+1,\ldots,b_1+p_1-1,\ldots,b_l,b_l+1,\ldots,b_l+p_l-1)=(\beta_1,\beta_2,\ldots,\beta_p)$ are distinct.
It is straightforward to verify the partial fraction decomposition
$$
\prod_{j=1}^{p}\frac{1}{\beta_j+x}=\frac{1}{(\betta+x)_1}=\sum\limits_{q=1}^{p}\frac{1}{B_q(\beta_q+x)},~~\text{where}~B_q=\prod\limits_{\substack{v=1\\v\ne q}}^{p}(\beta_v-\beta_q).
$$
Then
$$
\frac{(\b)_n}{(\b+\p)_n}=\frac{(\b)_\p}{(\b+n)_\p}=\frac{(\b)_{\p}}{(\betta+n)_1}=(\b)_\p\sum\limits_{q=1}^{p}\frac{(\beta_q)_n}{\beta_qB_q(\beta_q+1)_n}.
$$
Applying the definition of the generalized hypergeometric function and Theorem~\ref{thm:Fpred}, we obtain
\begin{multline}\label{eq:vector1}
(1-x)^a{}_{r+p+2}F_{r+p+1}\left.\!\!\left(\!\begin{matrix}a, \b,\f+\m\\ \b+\p,\f\end{matrix}\right\vert x\right)
=(1-x)^a\sum\limits_{q=1}^p
\frac{(\b)_{\p}}{\beta_qB_q}{}_{r+2}F_{r+1}\left.\!\!\left(\!\begin{matrix}a, \beta_q,\f+\m\\ \beta_q+1,\f\end{matrix}\right\vert x\right)
\\
=\frac{(\b)_{\p}}{(\f)_\m}\sum\limits_{q=1}^p\frac{(\f-\beta_q)_\m}{\beta_qB_q}{}_{2}F_{1}\left.\!\!\left(\!\begin{matrix}a, 1\\ \beta_q+1\end{matrix}\right\vert \frac{x}{x-1}\right)+(\b)_{\p}\sum\limits_{q=1}^p \frac{1}{\beta_qB_q}\sum\limits_{k=0}^{m-1}\delta_{kq}\sum\limits_{l=0}^k\SS_{k}^{(l)}\frac{(a)_lx^l}{(1-x)^l}
\\
=\frac{(\b)_{\p}}{(\f)_\m}\sum\limits_{q=1}^p\frac{(\f-\beta_q)_\m}{\beta_qB_q}{}_{2}F_{1}\left.\!\!\left(\!\begin{matrix}a, 1\\\beta_q+1\end{matrix}\right\vert\frac{x}{x-1}\right)
+(\b)_{\p}\sum\limits_{q=1}^p \frac{1}{\beta_qB_q}\sum\limits_{l=0}^{m-1}Y_l(\beta_q,\f,\m)\frac{(a)_lx^l}{(1-x)^l},
\end{multline}
where $\delta_{kq}$ are the coefficients of the polynomial $W_{m-1}(\beta_q,\f,\m;x)=\sum\nolimits_{k=0}^{m-1}\delta_{kq}x^k$, $Y_l$ is defined in (\ref{eq:Bl}) and $\beta_q=b+q-1$ is $q$-th component of the vector $\betta$.  Similarly, applying the second transformation yields:
\begin{multline}\label{eq:Kar1111}
(1-x)^{a-1}{}_{r+p+2}F_{r+p+1}\left.\!\!\left(\!\begin{matrix}a, \b,\f+\m\\ \b+\p,\f\end{matrix}\right\vert x\right)
=(1-x)^{a-1}\sum\limits_{q=1}^p \frac{(\b)_{\p}}{\beta_qB_q}{}_{r+2}F_{r+1}\left.\!\!\left(\!\begin{matrix}a, \beta_q,\f+\m\\ \beta_q+1,\f\end{matrix}\right\vert x\right)
\\
=\frac{(\b)_{\p}}{(\f)_\m}\sum\limits_{q=1}^p \frac{(\f-\beta_q)_{\m}}{\beta_qB_q}
{}_{2}F_{1}\left.\!\!\left(\!\begin{matrix}1,\beta_q+1-a\\\beta_q+1\end{matrix}\right\vert x\right)
+\sum\limits_{q=1}^p \frac{(\b)_{\p}}{\beta_qB_q}\sum\limits_{l=0}^{m-1}Y_l(\beta_q,\f,\m)\frac{(a)_lx^l}{(1-x)^{l+1}}.
\end{multline}
In both formulas the sum of the Gauss functions ${}_{2}F_{1}$ does not seem to collapse into a single hypergeometric function.  However, it does happen when $\b$ only contains one component. We formulate this result in the form of the following theorem.

\begin{theorem} \label{thm:Fpred1}
Suppose $p\in\N$. Then the following identity hold true\emph{:}
\begin{multline}\label{eq:Kar11}
(1-x)^a{}_{r+2}F_{r+1}\left.\!\!\left(\begin{matrix}a, b,
\f+\m\\ b+p,\f\end{matrix}\right\vert x\right)=\frac{T_{p-1}(0)}{\Gamma(b)(\f)_\m}{}_{p+1}F_{p}\left.\!\!\left(\begin{matrix}a, 1,-\llambda+1
\\ b+p,-\llambda\end{matrix}\right\vert \frac{x}{x-1}\right)
\\
+\sum\limits_{q=1}^p\frac{(-1)^{q-1}(b)_{p}}{(b+q-1)(q-1)!(p-q)!} \sum\limits_{l=0}^{m-1}Y_l(b+q-1,\f,\m)\frac{(a)_lx^l}{(1-x)^l},
\end{multline}
where $Y_l$ is defined in \emph{(\ref{eq:Bl})}, $\llambda=(\lambda_1,\ldots,\lambda_{p-1})$ are the roots of the polynomial
\begin{equation}\label{eq:Tp}
T_{p-1}(z)=\sum\limits_{q=1}^p\frac{(-1)^{q-1}(\f-b-q+1)_\m\Gamma(b+q-1)}{(q-1)!(p-q)!}(b+q+z)_{p-q}
\end{equation}
of degree $p-1$.  Furthermore,
\begin{multline}\label{eq:Kar12}
(1-x)^{a-1}{}_{r+2}F_{r+1}\left.\!\!\left(\begin{matrix}a, b,
\f+\m\\ b+p,\f\end{matrix}\right\vert x\right)
\!=\!\frac{\Gamma(b-a+1)T^*_{p-1}(0)}{\Gamma(b)(\f)_{\m}}{}_{p+1}F_{p}\left.\!\!\left(\begin{matrix}1,b+1-a,-\llambda^*+1\\ b+p,-\llambda^*\end{matrix}\right\vert x\right)
\\
+\sum\limits_{q=1}^p\frac{(-1)^{q-1}(b)_{p}}{(b+q-1)(q-1)!(p-q)!} \sum\limits_{l=0}^{m-1}Y_l(b+q-1,\f,\m)\frac{(a)_lx^l}{(1-x)^{l+1}},
\end{multline}
where  $\llambda^*=(\lambda_1^*,\ldots,\lambda^*_{p-1})$  are the roots of the polynomial
$$
T^*_{p-1}(z)=\sum\limits_{q=1}^p \frac{(-1)^{q-1}(\f-b-q+1)_\m\Gamma(b+q-1)}{\Gamma(b+q-a)(q-1)!(p-q)!}(b+q+z)_{p-q}(b+1-a+z)_{q-1}
$$
of degree $p-1$.
\end{theorem}

\textbf{Remark}.  It is instructive to compare identities (\ref{eq:Kar11}) and (\ref{eq:Kar12}) with the degenerate Miller-Paris transformations derived in \cite[Theorems~1~and~3]{KPReluts2018}.  One important difference is that the above theorem holds for any $p\in\N$, while in \cite{KPReluts2018} $p$ is restricted to the set $\{1,\ldots,m\}$.  Nevertheless, with a little effort one can make sure that (\ref{eq:Kar11}) and \cite[(16)]{KPReluts2018} are related by a rather simple rearrangement and  the polynomial  $T_{p-1}$ is a constant multiple of the polynomial $R_p$ from  \cite{KPReluts2018}. The same, however, does not hold for (\ref{eq:Kar12}), and the polynomial $T^*_{p-1}$ is not a constant multiple of $\hat{R}_p$ from \cite{KPReluts2018}.

\textbf{Proof.} If $\b=(b)$ we can represent the sum of hypergeometric functions in (\ref{eq:vector1}) as a single hypergeometric function of a higher order as follows. Using the definition of the hypergeometric function
\begin{equation}\label{eq:Kar113}
\sum\limits_{q=1}^p\frac{(\f-\beta_q)_\m}{\beta_qB_q}{}_{2}F_{1}\left.\!\!\!\left(\begin{matrix}a, 1\\\beta_q+1\end{matrix}\right\vert\:t\right)
=\sum\limits_{n=0}^\infty(a)_nt^n\sum\limits_{q=1}^p \frac{(\f-\beta_q)_\m}{\beta_qB_q(b+q)_n}.
\end{equation}
Applying the formula
\begin{equation}\label{eq:rr}
(b+q)_n=\frac{\Gamma(b+p)(b+p)_n}{\Gamma(b+q)(b+n+q)_{p-q}},
\end{equation}
we get
\begin{multline}\label{eq:Kar115}
\sum\limits_{n=0}^\infty(a)_nt^n\sum\limits_{q=1}^p \frac{(\f-\beta_q)_\m}{\beta_qB_q(b+q)_n}
=\frac{1}{\Gamma(b+p)}\sum\limits_{n=0}^\infty \frac{(a)_nt^n}{(b+p)_n}\sum\limits_{q=1}^p \frac{(\f-\beta_q)_\m}{\beta_qB_q}\Gamma(b+q)(b+n+q)_{p-q}
\\
=\frac{1}{\Gamma(b+p)}\sum\limits_{n=0}^\infty \frac{(a)_nt^n}{(b+p)_n}T_{p-1}(n),
\end{multline}
where $T_{p-1}(n)$ is the polynomial of degree $p-1$ defined in (\ref{eq:Tp}) in view of $\betta=(b,\ldots,b+p-1)$, $\beta_q=b+q-1$ and $B_q=(-1)^{q-1}(q-1)!(p-q)!$. Setting $\llambda=(\lambda_1,\ldots,\lambda_{p-1})$ to be the roots
of this polynomial, we can write
$$
T_{p-1}(n)=\frac{\Gamma(b+1)}{b(p-1)!}(\f-b)_\m (n-\llambda)_1=\frac{\Gamma(b)}{(p-1)!}(\f-b)_\m(-\llambda)_1\frac{(-\llambda+1)_n}{(-\llambda)_n}.
$$
Hence, it follows from (\ref{eq:Kar115}) that
\begin{multline}\label{eq:Kar116}
\sum\limits_{n=0}^\infty (a)_nt^n\sum\limits_{q=1}^p \frac{(\f-b-q+1)_\m}{\beta_qB_q(b+q)_n}
=\frac{\Gamma(b)}{\Gamma(b+p)}\sum\limits_{n=0}^\infty \frac{(a)_n t^n}{(b+p)_n}\frac{(\f-b)_\m(-\llambda)_1(-\llambda+1)_n}{(-\llambda)_n(p-1)!}
\\
=\frac{(\f-b)_\m(-\llambda)_1}{(b)_p(p-1)!}
{}_{p+1}F_{p}\left.\!\!\left(\begin{matrix}a, 1,-\llambda+1\\b+p,-\llambda\end{matrix}\right\vert t\right).
\end{multline}
Substituting this result into (\ref{eq:vector1}) with $t=x/(x-1)$ yields (\ref{eq:Kar11}).

The proof of the second transformation is similar.  We transform the first term in (\ref{eq:Kar1111}) using formula (\ref{eq:rr}) and the identity
$$
(b+n+q-a)_{p-q}=\frac{(b-a+n+1)_{p-1}}{(b-a+n+1)_{q-1}}
$$
as follows (keeping in mind that $\beta_q=b+q-1$, $B_q=(-1)^{q-1}(q-1)!(p-q)!$):
\begin{multline}\label{eq:Kar11114}
\frac{(b)_{p}}{(\f)_\m}\sum\limits_{q=1}^p \frac{(\f-\beta_q)_\m}{\beta_qB_q}{}_{2}F_{1}\left.\!\!\left(\begin{matrix}\beta_q-a+1,1\\\beta_q+1\end{matrix}\right\vert x\right)
=\frac{(b)_{\p}}{(\f)_\m}\sum\limits_{n=0}^\infty x^n\sum\limits_{q=1}^p \frac{(\f-\beta_q)_\m}{\beta_qB_q}\frac{(b+q-a)_n}{(b+q)_n}
\\
=\frac{(b)_{p}}{(\f)_\m}\sum\limits_{n=0}^\infty x^n\sum\limits_{q=1}^p \frac{(\f-\beta_q)_\m}{\beta_qB_q}\frac{\Gamma(b+q)(b+n+q)_{p-q}\Gamma(b+p-a)(b+p-a)_n}{\Gamma(b+p)(b+p)_n\Gamma(b+q-a)(b+n+q-a)_{p-q}}
\\
=\frac{\Gamma(b+p-a)}{\Gamma(b)(\f)_\m}\sum\limits_{n=0}^\infty\frac{(b+p-a)_nx^n}{(b+p)_n} \sum\limits_{q=1}^p \frac{(\f-\beta_q)_\m}{\beta_qB_q}\frac{\Gamma(b+q)(b+n+q)_{p-q}}{\Gamma(b+q-a)(b+n+q-a)_{p-q}}
\\
=\frac{\Gamma(b+p-a)}{\Gamma(b)(\f)_\m}\sum\limits_{n=0}^\infty\frac{(b+p-a)_nx^n}{(b+p)_n(b+n+1-a)_{p-1}}T^*_{p-1}(n)
\\
=\frac{1}{\Gamma(b)(\f)_\m}\sum\limits_{n=0}^\infty\frac{\Gamma(b+n-a+1)x^n}{(b+p)_n}T^*_{p-1}(n)
=\frac{\Gamma(b-a+1)}{\Gamma(b)(\f)_\m}\sum\limits_{n=0}^\infty\frac{(b-a+1)_nx^n}{(b+p)_n}T^*_{p-1}(n),
\end{multline}
where
$$
T^*_{p-1}(z)=\sum\limits_{q=1}^p \frac{(-1)^{q-1}\Gamma(b+q-1)(\f-b-q+1)_\m}{\Gamma(b+q-a)(q-1)!(p-q)!}
(b+z+q)_{p-q}(b+z+1-a)_{q-1}
$$
is a polynomial of degree $p-1$. Setting $\llambda^*=(\lambda_1^*,\ldots,\lambda^*_{p-1})$ to be the roots of this polynomial, we have
$$
T^*_{p-1}(n)=(-\llambda^*)_1\frac{(-\llambda^*+1)_n}{(-\llambda^*)_n}\sum\limits_{q=1}^p (-1)^{q-1}\frac{\Gamma(b+q-1)(\f-b-q+1)_{\m}}{\Gamma(b+q-a)(q-1)!(p-q)!}.
$$
Substituting this expression into (\ref{eq:Kar11114}), we get (\ref{eq:Kar12}).$\hfill\square$

The remark made  after  Theorem~\ref{thm:Fpred} implies that $p=2$ case of \eqref{eq:Kar11} is the same (modulo some rearrangement) as \cite[Corollary~5]{KPReluts2018}.
Setting $p=2$ in \eqref{eq:Kar12} we obtain the following

\begin{corollary}
Suppose $(b+1)(\f-b-1+\m)\ne b(\f-b-1)$ and  $(b-a+1)(\f-b-1+\m)\ne b(\f-b-1)$. Then the following identity holds:
\begin{multline*}
(1-x)^{a-1}{}_{r+2}F_{r+1}\left.\!\!\left(\begin{matrix}a, b,\f+\m\\ b+2,\f\end{matrix}\right\vert x\right)
\!=\!\frac{b\left[(\f-b)_{\m}-(\f-b-1)_{\m}\right]}{(\f)_{\m}}
{}_{3}F_{2}\left.\!\!\left(\begin{matrix}1,b+1-a,\lambda^*+1\\ b+2,\lambda^*\end{matrix}\right\vert x\right)
\\
+\sum\limits_{l=0}^{m-1}[(b+1)Y_l(b,\f,\m)-bY_l(b+1,\f,\m)]\frac{(a)_lx^l}{(1-x)^{l+1}},
\end{multline*}
where
$$
\lambda^*=(b-a+1)\frac{(b+1)(\f-b-1+\m)-b(\f-b-1)}{(b-a+1)(\f-b-1+\m)-b(\f-b-1)}
$$
and $Y_l$ is defined in \emph{(\ref{eq:Bl})}.
\end{corollary}

The following theorem extends Theorem~\ref{thm:Fpred} in a different direction: we add two free parameters to ${}_{r+2}F_{r+1}$ on the left hand side.
\begin{theorem} \label{thm:Trans3}
Suppose $(e-d-m+1)_{m-1}\ne0$. Then following identity holds:
\begin{multline*}
{}_{r+3}F_{r+2}\left.\!\!\left(\!\begin{matrix}a,d,b,\f+\m\\e, b+1,\f\end{matrix}\right\vert x\right)
=\frac{(\f-b)_\m}{(\f)_\m}{}_{3}F_{2}\left.\!\!\left(\!\begin{matrix}a,d,b\\ e,b+1\end{matrix}\right\vert x\right)
\\
+\frac{(\f)_{\m}-(\f-b)_\m}{(\f)_\m}(1-x)^{-a}{}_{m+1}F_{m}\left.\!\!\left(\!\begin{matrix}a,e-d-m+1,\llambda+1\
\\e,\llambda\end{matrix}\right\vert \frac{x}{x-1}\right),
\end{multline*}
where $\llambda$ is the vector of zeros of the polynomial
\begin{equation}\label{eq:Lm}
L_{m-1}(t)=L_{m-1}(e,d,b,c,\f,\m;t)=\sum\limits_{k=0}^{m-1}(d)_kY_{k}(b,\f,\m)(t)_{k}(e-d-m+1-t)_{m-1-k},
\end{equation}
 and $Y_k(b,\f,\m)$ is given in \emph{(\ref{eq:Bl})}.
If, in addition $(e-a-m+1)_{m-1}\ne0$ and $(1+a+d-e)_{m-1}\ne0$, then
\begin{multline*}
{}_{r+3}F_{r+2}\left.\!\!\left(\!\begin{matrix}a,d,b,\f+\m\\e, b+1,\f\end{matrix}\right\vert x\right)
=\frac{(\f-b)_\m}{(\f)_\m}{}_{3}F_{2}\left.\!\!\left(\!\begin{matrix}a,d,b\\e,b+1\end{matrix}\right\vert x\right)
\\
+\frac{(\f)_{\m}-(\f-b)_\m}{(\f)_\m}(1-x)^{e-a-d-m+1}
{}_{m+1}F_{m}\left.\!\!\left(\!\begin{matrix}e-a-m+1,e-d-m+1,\llambda^*+1\\e,\llambda^*\end{matrix}\right\vert x\right),
\end{multline*}
where $\llambda^*$ is the vector of zeros of the polynomial
\begin{equation}\label{eq:Lmhat}
\hat{L}_{m-1}(t)=\sum\limits_{k=0}^{m-1}\frac{(-1)^kY_k(b,\f,\m)(a)_k(d)_k(t)_k}{(e-a-m+1)_k(e-d-m+1)_k}
{}_{3}F_{2}\!\left(\begin{matrix}-m+1+k,t+k,e-a-d-m+1\\e-a-m+1+k,e-d-m+1+k\end{matrix}\right).
\end{equation}
\end{theorem}
\textbf{Proof.} Let $\ggamma=(\gamma_1,\ldots,\gamma_{m-1})$ be the roots of the polynomial $W_{m-1}(x)=W_{m-1}(b,\f,\m;x)$ defined in (\ref{eq:pol}). Its definition implies that the leading coefficient of $W_{m-1}(x)$ equals $b/(\f)_m$, while the free term is given by $W_{m-1}(0)=((\f)_{\m}-(\f-b)_{\m})/(\f)_{\m}$ . Hence,
\begin{multline}\label{eq:Wnfactorization}
W_{m-1}(n)=\frac{b(-\ggamma)_1(n-\ggamma)_1}{(\f)_m (-\ggamma)_1}=\frac{b (-\ggamma)_1}{(\f)_m}\frac{(-\ggamma+1)_n}{(-\ggamma)_n}
\\
=W_{m-1}(0)\frac{(-\ggamma+1)_n}{(-\ggamma)_n}
=\frac{((\f)_{\m}-(\f-b)_{\m})(-\ggamma+1)_n}{(\f)_{\m}(-\ggamma)_n}.
\end{multline}
By definition of the generalized hypergeometric function this leads to 
\begin{multline*}
{}_{r+3}F_{r+2}\left.\!\!\left(\!\begin{matrix}a,d, b,\f+\m\\ e,b+1,\f\end{matrix}\right\vert x\right)
=\frac{(\f-b)_\m}{(\f)_\m}\sum\limits_{n=0}^\infty \frac{(a)_n (d)_n (b)_n x^n}{n!(e)_n(b+1)_n}
+\sum\limits_{n=0}^\infty \frac{(a)_n(d)_nx^n W_{m-1}(n)}{n!(e)_n}
\\
=\frac{(\f-b)_\m}{(\f)_\m}{}_{3}F_{2}\left.\!\!\left(\!\begin{matrix}a,d,b\\e,b+1\end{matrix}\right\vert x\right)
+\sum\limits_{n=0}^\infty \frac{(a)_n  (d)_n x^n W_{m-1}(n)}{n! (e)_n}
\\
=\frac{(\f-b)_\m}{(\f)_\m}{}_{3}F_{2}\left.\!\!\left(\!\begin{matrix}a,d,b\\e,b+1\end{matrix}\right\vert x\right)
+\frac{(\f)_m-(\f-b)_\m}{(\f)_\m}{}_{m+1}F_{m}\left.\!\!\left(\!\begin{matrix}a,d,-\ggamma+1\\e,-\ggamma\end{matrix}\right\vert x\right).
\end{multline*}
It remains to apply the Miller-Paris transformations (\ref{eq:KRPTh1-1}) and (\ref{eq:KRPTh1-2}) (or (\ref{eq:MP1}) and (\ref{eq:MP2})) to the function
\begin{equation*}
{}_{m+1}F_{m}\left.\!\!\left(\!\begin{matrix}a,d,-\ggamma+1\\e,-\ggamma\end{matrix}\right\vert x\right).
\end{equation*}
Note the change of notation $b\to d,$ $c\to e$, $m\to m-1$ as compared to (\ref{eq:KRPTh1-1}), (\ref{eq:KRPTh1-2}). To give explicit  
formulas for the characteristic polynomials  (\ref{eq:Qm}) and (\ref{eq:Qmhat}) use (\ref{eq:pol}) and (\ref{eq:Wnfactorization}) to get
\begin{multline*}
{}_{r+3}F_{r+2}\!\!\left(\!\begin{matrix}-k,b,\f+\m\\ b+1,\f\end{matrix}\right)=
\frac{(\f-b)_{\m}}{(\f)_\m}{}_{2}F_{1}\!\!\left(\begin{matrix}-k,b\\ b+1\end{matrix}\right)
+\frac{(\f)_{\m}-(\f-b)_\m}{(\f)_\m}{}_{m+1}F_{m}\!\!\left(\!\begin{matrix}-k,-\ggamma+1\\-\ggamma\end{matrix}\right)
\\
=\frac{(\f-b)_{\m}}{(\f)_{\m}}\frac{k!}{(b+1)_k}+\frac{(\f)_{\m}-(\f-b)_{\m}}{(\f)_{\m}}{}_{m+1}F_{m}\!\!\left(\!\begin{matrix}-k,-\ggamma+1\\\-\ggamma\end{matrix}\right),
\end{multline*}
where the Chu-Vandermonde identity was applied in the second equality. Comparing this formula with (\ref{eq:Ckr}) and (\ref{eq:Bl}) we immediately see that
$$
C_{k,r}(-\ggamma,\mathbf{1})=\frac{(-1)^k}{k!}{}_{m+1}F_{m}\left.\!\!\left(\!\begin{matrix}-k,-\ggamma+1\\-\ggamma\end{matrix}\right\vert x\right)=\frac{(\f)_{\m}}{(\f)_{\m}-(\f-b)_{\m}}Y_k(b,\f,\m).
$$
Substituting this expression into (\ref{eq:Qm}) and (\ref{eq:Qmhat}) and canceling constant factors we arrive at (\ref{eq:Lm}) and (\ref{eq:Lmhat}), respectively.
$\hfill\square$

Taking $r=1$, $m=2$ in Theorem~\ref{thm:Trans3} after some elementary computations we arrive at
\begin{corollary} \label{cr:f2b1}
Suppose $e-d-1\ne0$. Then following identities hold:
\begin{multline*}
f(f+1){}_{4}F_{3}\left.\!\!\left(\!\begin{matrix}a,d,b,f+2\\e, b+1,f\end{matrix}\right\vert x\right)
-(f-b)(f-b+1){}_{3}F_{2}\left.\!\!\left(\!\begin{matrix}a,d,b\\ e,b+1\end{matrix}\right\vert x\right)
\\
=b(2f-b+1)(1-x)^{-a}{}_{3}F_{2}\left.\!\!\left(\!\begin{matrix}a,e-d-1,\lambda+1\
\\e,\lambda\end{matrix}\right\vert \frac{x}{x-1}\right)
\\
=b(2f-b+1)(1-x)^{e-a-d-1}
{}_{3}F_{2}\left.\!\!\left(\!\begin{matrix}e-a-1,e-d-1,\lambda^*+1\\e,\lambda^*\end{matrix}\right\vert x\right),
\end{multline*}
where
$$
\lambda=\frac{(2f-b+1)(e-d-1)}{2f-b-d+1},~~~
\lambda^*=\frac{(2f-b+1)(e-a-1)(e-d-1)}{ad+(2f-b+1)(e-a-d-1)}.
$$
For the second equality the additional restrictions $e-a-1\ne0$ and $1+a+d-e\ne0$ must be imposed.
\end{corollary}

Similarly, taking $r=2$, $m_1=m_2=1$ in Theorem~\ref{thm:Trans3} we get
\begin{corollary} \label{cr:f1f1b1}
Suppose $e-d-1\ne0$. Then following identities hold:
\begin{multline*}
(f_1f_2){}_{5}F_{4}\left.\!\!\left(\!\begin{matrix}a,d,b,f_1+1,f_2+1\\e, b+1,f_1,f_2\end{matrix}\right\vert x\right)
-(f_1-b)(f_2-b){}_{3}F_{2}\left.\!\!\left(\!\begin{matrix}a,d,b\\ e,b+1\end{matrix}\right\vert x\right)
\\
=b(f_1+f_2-b)(1-x)^{-a}{}_{3}F_{2}\left.\!\!\left(\!\begin{matrix}a,e-d-1,\lambda+1\
\\e,\lambda\end{matrix}\right\vert \frac{x}{x-1}\right)
\\
=b(f_1+f_2-b)(1-x)^{e-a-d-1}
{}_{3}F_{2}\left.\!\!\left(\!\begin{matrix}e-a-1,e-d-1,\lambda^*+1\\e,\lambda^*\end{matrix}\right\vert x\right),
\end{multline*}
where
$$
\lambda=\frac{(f_1+f_2-b)(e-d-1)}{f_1+f_2-b-d},~~~
\lambda^*=\frac{(f_1+f_2-b)(e-a-1)(e-d-1)}{ad+(f_1+f_2-b)(e-a-d-1)}.
$$
For the second equality the additional restrictions $e-a-1\ne0$ and $1+a+d-e\ne0$ must be imposed.
\end{corollary}


\begin{thebibliography}{99}
\bibitem{AAR} G.E.\:Andrews, R.\:Askey and R.\:Roy, Special functions, Cambridge University Press, 1999.

\bibitem{BealsWong} R.\:Beals and R.\:Wong, Special Functions and Orthogonal Polynomials,
Cambridge Studies in Advanced Mathematics   (No. 153), Cambridge University Press, 2016.

\bibitem{Gasper}G.\:Gasper, Summation formulas for basic hypergeometric series. SIAM J. Math. Anal. 12(1981), 196--200.

\bibitem{Chu1}W.\:Chu, Partial fractions and bilateral summations. Journal of Mathematical Physics vol.35(1994), 2036.

\bibitem{Chu2}W.\:Chu, Erratum: Partial fractions and bilateral summations. Journal of Mathematical Physics. vol.36 (1995), 5198.

\bibitem{Karlson}P.W.\:Karlsson, Hypergeometric functions with integral parameter differences. J. Math. Phys. 12(1971), 270--271.

\bibitem{KarpJMS2015}D.\:Karp, Representations and inequalities for generalized hypergeometric functions, Journal of Mathematical Sciences, Volume 207, Issue 6(2015), 885--897.

\bibitem{KLJAT2017}D.\:Karp and J.L.\:L\'{o}pez, Representations of hypergeometric functions for arbitrary values of the parameters and their use,
Journal of Approximation Theory, Volume 218(2017), 42--70.

\bibitem{KMT2018} D.B.\:Karp, Yu.B.\:Melnikov and I.V.\:Turuntaeva, On the properties of special functions generating the kernels of certain integral operators,
Journal of Computaional and Applied Mathematics, submitted, 2018. Preprint, arXiv:1804.03982.

\bibitem{KPSIGMA} D.\:Karp and E.\:Prilepkina, Hypergeometric differential equation and new identities for the coefficients of
N{\o}rlund and B\"{u}hring, SIGMA 12 (2016), 052, 23 pages.

\bibitem{KPITSF2017} D.\:Karp and E.\:Prilepkina,  Applications of the Stieltjes and Laplace transform representations of the hypergeometric functions, 
Integral Transforms and Special Functions, volume 28, no.10 (2017), 710--731.

\bibitem{KPReluts2018}D.B.\:Karp and E.G.\:Prilepkina, Degenerate Miller-Paris transformations, submitted to Results in Mathematics, 2018. Preprint arXiv:1806.00208

\bibitem{KPITSF2018}D.B.\:Karp and E.G.\:Prilepkina, Extensions of Karlsson–Minton summation theorem and some consequences of the first Miller–Paris transformation,
Integral Transforms and Special Functions, Vol. 29, Issue 12 (2018), 955--970.


\bibitem{KRP2014}Y.S.\:Kim, A.K.\:Rathie and R.B.\:Paris, On two Thomae-type transformations for hypergeometric series with integral parameter differences,
Math. Commun. 19(2014), 111--118.


\bibitem{LVW}J.\:Letessier, G.\:Valent  and J.\:Wimp, Some differential equations satisfied by hypergeometric functions. Approximation and Computation (International Series of Numerical Mathematics vol 119). Cambridge, MA: Birkh\"{a}user, 1994, 371--81.

\bibitem{LukeBook} Y.L.\:Luke, The special functions and their approximations. Volume 1. Academic Press, 1969.

\bibitem{Miller2005}A.R.\:Miller, A summation formula for Clausen's series ${}_3F_2(1)$ with an application to Goursat’s function ${}_2F_2(x)$,
J. Phys. A: Math. Gen. 38 (2005), 3541--3545.

\bibitem{MP2011}A.R.\:Miller and R.B.\:Paris, Euler-type transformations for the generalized hypergeometric function ${}_{r+2}F_{r+1}(x)$,
Z. Angew. Math. Phys., Volume 62, Issue 1(2011),  31--45.

\bibitem{MP2012R} A.R.\:Miller, R.B.\:Paris, On a result related to transformations and summations of generalized hypergeometric series, 
Mathematical communications, Volume 17, 2012, 205--210.

\bibitem{MP2013}A.R.\:Miller and R.B.\:Paris, Transformation formulas for the generalized hypergeometric function with integral parameter differences, 
Rocky Mountain Journal Of Mathematics Volume 43, Number 1 (2013), 291--327.

\bibitem{MS2010}A.R.\:Miller and  H.M.\:Srivastava, Karlsson-Minton summation theorems for the generalized hypergeometric series of unit argument, 
Integral Transfroms and Special Functions, Vol. 21, No. 8(2010), 603--612.

\bibitem{Minton}Minton\:B.M. Generalized hypergeometric functions at unit argument. J. Math. Phys. 12(1970), 1375--1376.

\bibitem{Norlund} N.E.\,N{\o}rlund, Hypergeometric functions, Acta Mathematica, volume 94(1955), 289--349.

\bibitem{NIST}F.W.J.\:Olver, D.W.\:Lozier,  R.F.\:Boisvert and C.W.\:Clark (Eds.) NIST Handbook of Mathematical Functions, Cambridge
University Press, 2010.

\bibitem{PBM3} A.P.\:Prudnikov, Yu.A.\:Brychkov, O.I.\:Marichev, Integrals and Series: More Special Functions,  Volume~3, Gordon and Breach Science Publishers, 
New York, 1990.

\bibitem{RJRJR} K.S.\:Rao, J.\:Van der Jeugt, J.\:Raynal, R.\:Jagannathan, V.\:Rajeswari, Group theoretical basis for the terminating $3F2(1)$ series,
J. Phys. A Math. Gen. {\bf 25} (1992), 861--876.

\bibitem{Rosengren1}H.\:Rosengren, Reduction Formulas for Karlsson–Minton-Type Hypergeometric Functions. Constr. Approx. 20(2004), 525--548.

\bibitem{Rosengren2}H.\:Rosengren,  Karlsson–Minton type hypergeometric functions on the root system $C_n$. J. Math. Anal. Appl. 281 (2003), 332--345.

\bibitem{Schlosser1}M.\:Schlosser, Multilateral transformations of $q$-series with quotients of parameters that are nonnegative integral powers of $q$, 
in: $q$-Series with Applications to Combinatorics, Number Theory, and Physics (B. C. Berndt and K. Ono, eds.). Amer. Math. Soc. Contemp. Math. 291(2001), 203--227.

\bibitem{Schlosser2}M.\:Schlosser, Elementary derivations of identities for bilateral basic hypergeometric series. Selecta Math. (N.S.). 9:1(2003), 119--159.

\bibitem{Seaborn}J.B.\:Seaborn, Hypergeometric Functions and Their Applications, Springer, 1991.
\end{thebibliography}
\end{document}